\documentclass[11pt, a4paper]{amsart}

\title{On the Topology of Fano Smoothings}
\author{Tom Coates}
\author{Alessio Corti}
\author{Genival da Silva Jr}
\usepackage{graphicx}
\usepackage{hyperref}

\usepackage{amsmath, amsthm, amssymb}
\usepackage[all]{xy}
\usepackage{cite}
\usepackage[margin=3cm]{geometry}

\newtheorem{theorem}{Theorem}

\newtheorem{definition}{Definition}

\newcommand{\ZZ}{\mathbb{Z}}
\newcommand{\QQ}{\mathbb{Q}}
\newcommand{\RR}{\mathbb{R}}
\newcommand{\CC}{\mathbb{C}}
\newcommand{\PP}{\mathbb{P}}

\newcommand{\cX}{\mathcal{X}}
\newcommand{\cY}{\mathcal{Y}}
\newcommand{\cO}{\mathcal{O}}

\DeclareMathOperator{\GL}{GL}

\DeclareMathOperator{\Pic}{Pic}
\DeclareMathOperator{\Proj}{Proj}
\DeclareMathOperator{\rk}{rk}

\usepackage{tikz}
\usetikzlibrary{positioning}
\usepackage{tikz-3dplot}

\begin{document}

\begin{abstract}
	Suppose that $X$ is a Fano manifold that corresponds under Mirror Symmetry to a Laurent polynomial $f$, and that $P$ is the Newton polytope of $f$.  In this setting it is expected that there is a family of algebraic varieties over the unit disc with general fiber $X$ and special fiber the toric variety defined by the spanning fan of $P$.  Building on recent work and conjectures by Corti--Hacking--Petracci, who construct such families of varieties, we determine the topology of the general fiber from combinatorial data on~$P$.  This provides evidence for the Corti--Hacking--Petracci conjectures, and verifies that their construction is compatible with expectations from Mirror Symmetry.
\end{abstract}

\maketitle

\section{Introduction}

There has been recent interest in the classification of Fano manifolds via Mirror Symmetry~\cite{CoatesCortiGalkinGolyshevKasprzyk2013,CoatesCortiGalkinKasprzyk2016,AkhtarCoatesCortiHeubergerKasprzykOnetoPetracciPrinceTveiten2016}.  For us, an $n$-dimensional Fano manifold $X$ corresponds under Mirror Symmetry to a Laurent polynomial $f \in \CC[x_1^{\pm 1},\ldots,x_n^{\pm 1}]$ if the regularized quantum period of $X$, which is a generating function for certain genus-zero Gromov--Witten invariants of $X$, coincides with classical period $\pi_f$ of $f$:
\[
\pi_f(t) = \frac{1}{(2 \pi \mathrm{i})^n} \int_{S^1 \times \cdots \times S^1} \frac{1}{1 - t f} \frac{dx_1}{x_1} \cdots \frac{dx_n}{x_n}.
\]
If a Fano manifold $X$ corresponds under Mirror Symmetry to a Laurent polynomial $f$ then it is expected that there is a degeneration $\cX \to \Delta$, where $\Delta \subset \CC$ is a disc, such that the general fiber is $X$ and the special fiber is the toric variety $X_f$ defined by the spanning fan of the Newton polytope of $f$.  In general, $X_f$ will be singular.  It is natural to ask whether one can determine the Betti numbers of the Fano manifold $X$ from its mirror partner~$f$.  We will show that, in dimension~$3$, the answer to this question is `yes'.

The key new ingredient presented here is work by Corti--Hacking--Petracci~\cite{CortiHackingPetracci2019}, in preparation and still partly conjectural, which goes a long way towards establishing the expected picture described above in dimension~$3$.  Corti--Hacking--Petracci construct, given a 3-dimensional reflexive polytope $P$ decorated (as described below) with `decomposition data', a smoothing $\cX \to \Delta$ of the toric variety $X_P$ defined by the spanning fan of $P$.  The general fiber of this smoothing is a three-dimensional Fano manifold $X$.  The decomposition data also determine a Laurent polynomial $f$ with Newton polytope $P$, so that $X_f = X_P$.  In this paper we show that the Betti numbers of $X$ depend on the decomposition data only via $f$.  That is, for each of the (many) choices of decomposition data that give rise to the same Laurent polynomial $f$, the Betti numbers of the corresponding smoothing are the same.  Furthermore, these coincide with the Betti numbers of the Fano manifold that corresponds to $f$ under Mirror Symmetry.

We proceed as follows.  Corti--Hacking--Petracci construct, given a choice of decomposition data for a reflexive polytope $P$, a toric partial resolution $\pi \colon Y \to X_P$ and a family $\cY \to \Delta$ with special fiber $Y$ and fiber  over a general point $t \in \Delta$ a weak Fano manifold~$Y_t$.  Contracting finitely many $(-1,-1)$ curves in $Y_t$ gives a resolution $\pi_t \colon Y_t \to X_t$, where $X_t$ is a Fano variety with ordinary double points, and this fits into a diagram
\begin{equation}
    \label{eq:main_diagram}
    \begin{aligned}
        \xymatrix{
            & & Y_t \ar@{|~}[rr] \ar[d]_{\pi_t} & & Y \ar[d]^\pi \\
            X_\eta \ar@{|~}[rr] &&  X_t \ar@{|~}[rr] & & X_P  
        }
    \end{aligned}
\end{equation}    
where the arrow $\xymatrix{A \ar@{|~}[r] & B}$ means that $A$ is the general fiber in a family over $\Delta$ with special fiber $B$.  Here $\xymatrix{X_\eta \ar@{|~}[r] &  X_t}$ is Namikawa's smoothing of Fano varieties with ordinary double points.  The Fano variety $X_\eta$ is our desired smoothing of $X_P$, and the diagram above allows us to compute its Betti numbers.  The central fiber $Y$ is a toric variety, so we know its cohomology groups explicitly; we can compute the Betti numbers of $Y_t$ by analysing the vanishing cycles of the degeneration $\xymatrix{Y_t \ar@{|~}[r] & Y}$\! and, since the left-hand part of the diagram is a conifold transition from $Y_t$ to $X_\eta$, this determines the Betti numbers of $X_\eta$.

We begin by reviewing the cohomology of toric varieties and the vanishing cycle exact sequence.  We then explain in \S\ref{sec:smoothing} how to compute the Betti numbers of the smoothing $X$ from the decomposition data, and give examples in \S\ref{sec:examples}.  In \S\ref{sec:main_result} we prove that the Betti numbers of $X$ depend on the decomposition data only via the Laurent polynomial~$f$ determined by those data, and in \S\ref{sec:computations} we verify that these Betti numbers of $X$ coincide with the Betti numbers of the Fano manifold that corresponds to $f$ under Mirror Symmetry.

\section{Cohomology and Vanishing Cycles}

\subsection{The cohomology of toric varieties}

We will compute the Betti numbers of the fiber $Y$ in diagram \eqref{eq:main_diagram} using the fact that it is a toric variety\footnote{We learned the statement of Theorem~\ref{thm:toric_Betti} from Andrea Petracci.}. 

\begin{theorem}[cf.~\cite{Jordan1997}, Proposition~3.5.3] \label{thm:toric_Betti}
	Let $\Sigma$ be a complete fan in a three-dimensional lattice and let $X_\Sigma$ be the toric threefold defined by $\Sigma$. Let $d_i$ denote the number of $i$-dimensional cones in $\Sigma$, and let $b_i$ denote the $i$th Betti number of $X_\Sigma$.  Then:
	\begin{align*}
		& d_1 - d_2 + d_3 = 2 && b_2 = \rk \Pic(X_\Sigma) \\
		& b_0 = b_6 = 1 && b_3 = \rk \Pic(X_\Sigma) - d_2 + 2d_1 - 3 \\
		& b_1 = b_5 = 0 && b_4 = d_1 - 3
	\end{align*}
	and the Euler characteristic of $X$ is $d_3$.
 \end{theorem}
 
\subsection{The vanishing cycle exact sequence}

We will compute the Betti numbers of the fiber $Y_t$ in diagram \eqref{eq:main_diagram} by analysing the vanishing cycles for the degeneration $\xymatrix{Y_t \ar@{|~}[r] & Y}$. Consider a complex variety $\cY$, a disc $\Delta \subset \CC$, and a projective morphism $f:\cY\rightarrow \Delta$.  Let $\Delta^*=\Delta \setminus \{0\}$ be the punctured disc, and  
\[
\{0\}\xrightarrow{i_0}\Delta\xleftarrow{j_0} \Delta^*
\]
be the natural inclusions. Denote the fiber over $t \in \Delta^*$ by $Y_t$, and the fiber over $0 \in \Delta$ by $Y$. Choose a universal covering map $p_0:\widetilde{\Delta^*}\rightarrow \Delta^*$, and consider the diagram
\[ 
\xymatrix{
Y\ar[rr]^i \ar[d] & & \cY \ar[d]^f & & \cY \setminus Y \ar[d] \ar[ll]_j & & \widetilde{\cY\setminus Y} \ar[d] \ar[ll]_p \\
\{0\} \ar[rr]^{i_0} & & \Delta & & \Delta^* \ar[ll]_{j_0} & & \widetilde{\Delta^*}. \ar[ll]_{p_0}
}
\]
Let $\mathcal{S}$ be a stratification for $\cY$ and suppose that $\textbf{F}^\bullet\in \textbf{D}^b_\mathcal{S}(\cY)$. The \emph{nearby sheaf} is defined~\cite[Expos\'e~I]{SGA7} to be the complex $$\psi_f \textbf{F}^\bullet = i^{*} \textbf{R}(j\circ p)_*(j\circ p)^* \textbf{F}^\bullet. $$
By adjunction, there is a natural map $$i^{*}\textbf{F}^\bullet\rightarrow \psi_f \textbf{F}^\bullet.$$
The \emph{sheaf of vanishing cycles} $\phi_f \textbf{F}^\bullet$ is the cone on this map (ibid.), and there is a distinguished triangle:
\begin{equation}\label{short_van}
   i^{*}\textbf{F}^\bullet\rightarrow \psi_f \textbf{F}^\bullet\rightarrow \phi_f \textbf{F}^\bullet\xrightarrow{+1}
\end{equation}

Consider now the cohomology sheaves $\textbf{H}^i(\psi_f \textbf{F}^\bullet)$ and $\textbf{H}^i(\phi_f \textbf{F}^\bullet)$ -- these are complexes of sheaves on $Y_0$ -- and their stalks at $y \in Y_0$.    Embed $\cY$ into an affine space and let $B(y,\epsilon)$ be the open ball of radius $\epsilon$ around $y$. Then for sufficiently small $\epsilon > 0$ and for all $t\in\Delta^*$ such that $|t|<\epsilon$, we have
\[
\textbf{H}^i(\psi_f \textbf{F}^\bullet)_y\cong \mathbb{H}^i(F_{f,y},\textbf{F}^\bullet)
\]
where $F_{f,y}$ is the Milnor fiber of $f$ at $y$, $F_{f,y}:=B(y,\epsilon)\cap \cY\cap f^{-1}(t)$.  Similarly,
\[
\textbf{H}^i(\phi_f \textbf{F}^\bullet)_y\cong \mathbb{H}^{i+1}(B(y,\epsilon)\cap\cY ,F_{f,y};\textbf{F}^\bullet).
\]
Now consider the constant sheaf $\mathbb{Q}_\cY$ as a complex concentrated in degree~$0$. Taking stalks of the hypercohomology of the distinguished triangle in \eqref{short_van}, we get
\begin{equation} \label{eq:van}
	\cdots \rightarrow
    H^i(Y,\mathbb{Q})\rightarrow
	H^i(Y_t,\mathbb{Q})\rightarrow
	H^i_v(Y_t,\mathbb{Q})\rightarrow
    H^{i+2}(Y,\mathbb{Q})\rightarrow 
	\cdots
\end{equation}
where $H^i_v(Y_t,\mathbb{Q})$ is the subspace in $H^i(Y_t,\mathbb{Q})$ generated by vanishing cycles, that is, cycles in the kernel of the natural map $H^i(Y_t,\mathbb{Q})\rightarrow H^i(Y_0,\mathbb{Q})$.  

\section{Smoothing Toric Fano Threefolds}
\label{sec:smoothing}

Every three-dimensional Gorenstein toric Fano variety is the toric variety $X_P$ defined by the spanning fan of a three-dimensional reflexive polytope $P$.  This gives a one-to-one correspondence between three-dimensional Gorenstein toric Fano varieties up to isomorphism and three-dimensional reflexive polytopes up to $\GL(3,\ZZ)$-equivalence.  In general such a toric variety $X_P$ is singular.  As mentioned in the Introduction, Corti--Hacking--Petracci construct, starting from a three-dimensional reflexive polytope $P$ decorated with some additional data, a smoothing $\cX \to \Delta$ of $X_P$.  In this section we describe their construction and the additional data required.

\begin{definition}
	Let $n$ be an integer such that $n \geq -1$.  An \emph{$A_n$~triangle} is a lattice polygon that is $\ZZ^2 \rtimes \GL(\ZZ^2)$-equivalent to the polygon in $\QQ^2$ with vertices $(0,0)$, $(0,1)$, and $(n+1,0)$
\end{definition}

\begin{paragraph}{\textbf{Example}}
	The standard two-dimensional simplex is an $A_0$-triangle.
\end{paragraph}

\medskip

\begin{paragraph}{\textbf{Example}}
	An $A_{-1}$-triangle is a line segment of unit length.
\end{paragraph}

\begin{definition} \label{def:admissible}
	Let $F$ be a lattice polygon.  An \emph{admissible Minkowski decomposition} of $F$ is a Minkowski decomposition $F = F_1 + \cdots + F_k$ of $F$ as a sum of lattice polygons $F_j$ such that each $F_j$ is an $A_n$-triangle for some $n \geq -1$.  (Here $n$ can depend on $j$.)  We consider admissible Minkowski decompositions of $F$ that differ by reordering and translation of the summands to be equivalent.
\end{definition}

We now introduce certain polyhedral decompositions of lattice polygons.  Suppose that $F = F_1 + \cdots + F_k$ is an admissible Minkowski decomposition of the lattice polygon $F$.  Recall that the \emph{Cayley polytope} $C_{F_1,\ldots,F_k}$ is the convex hull of
\[
F_1 + e_1, F_2 + e_2, \ldots, F_k + e_k
\]
in $L \oplus \ZZ e_1 \oplus \cdots \oplus \ZZ e_k$, where $L$ is the lattice containing $F$.  The map
\begin{equation}
	\label{eq:Cayley_slice}
	\begin{aligned}
		F & \longrightarrow C_{F_1,\ldots,F_k} \cap \left(L_{\RR} \times \left\{\left(\textstyle\frac{1}{k},\cdots,\frac{1}{k}\right)\right\}\right) \\
		v & \longmapsto \textstyle \frac{1}{k}(v+e_1 + \cdots + e_k)
	\end{aligned}
\end{equation}
is bijective, and this allows us to obtain polyhedral decompositions of $F$ from polyhedral decompositions of the Cayley polytope $C_{F_1,\ldots,F_k}$.

\begin{definition} \label{def:subdivision}
	Let $F = F_1 + \cdots + F_k$ be an admissible lattice Minkowski decomposition as above.  A \emph{regular fine mixed subdivision} of $F$ is a subdivision of $F$ induced, via \eqref{eq:Cayley_slice}, by a regular unimodular triangulation of the Cayley polytope $C_{F_1,\ldots,F_k}$.  We say that such a subdivision is \emph{subordinate to the Minkowski decomposition} $F = F_1 + \cdots + F_k$.
\end{definition}

\noindent There is a whole body of theory here, which we will not discuss,
concerning polyhedral subdivisions which may be neither fine nor mixed: an
introduction to these topics can be found in the extremely beautiful book by
De~Loera et al.~\cite{DeLoeraRambauSantos2010}. We will only consider regular
fine mixed subdivisions.

For the rest of this section, fix a three-dimensional reflexive polytope $P$.  Corti--Hacking--Petracci consider the polytope $P$ together with  \emph{decomposition data}.  This is, for each facet $F$ of $P$:\label{choices}
\begin{itemize}
\item[(A)] a choice of admissible Minkowski decomposition $F = F_1 + \cdots + F_k$;
\item[(B)] a choice of regular fine mixed subdivision of $F$ subordinate to (A);
\end{itemize}
satisfying a condition that we now describe.  Note that by taking cones over the regular fine mixed subdivisions of each facet, we obtain a complete fan $\Sigma$ that refines the spanning fan of $P$, and thus a toric crepant partial resolution $Y \to X_P$.   Recall first that irreducible toric curves in $X_P$ (respectively in~$Y$) correspond to two-dimensional cones in the spanning fan of~$P$ (respectively in the fan~$\Sigma$).  Thus irreducible toric curves in $X_P$ correspond to edges of $P$, and irreducible toric curves in $Y$ correspond to edges in the polyhedral subdivision of the boundary of $P$ determined by the decomposition data.

\begin{theorem}[Corti--Hacking--Petracci]
	The singularities of the toric partial resolution $Y$ of $X_P$ are at worst quasi-ordinary double points (qODPs).
\end{theorem}

This amounts to the statement that each polygon in the regular fine mixed subdivision (B) above is either a standard $2$-simplex or a quadrilateral with sides of unit length. The curves in the toric $1$-skeleton of $Y$, therefore, either meet at $3$-valent vertices -- which are the torus-fixed points on $Y$ corresponding to the cones over the $2$-simplices -- or at $4$-valent vertices, which are the torus-fixed points on $Y$ corresponding to the cones over the quadrilaterals. The $3$-valent vertices are smooth points on $Y$, and the $4$-valent vertices are the qODPs.

Let $\Gamma_e$ denote the irreducible toric curve in $X_P$ determined by the edge $e$ of $P$, and let $\widetilde{\Gamma}_e$ denote the set of irreducible toric
curves in $Y$ that map dominantly to $\Gamma_e$ under the partial resolution $Y \to X_P$. The fact that each polygon in the subdivision (B) above is either a standard $2$-simplex or a quadrilateral with sides of unit length implies that $|\widetilde{\Gamma}_e| = \ell(e)$, the lattice length of the edge $e$.

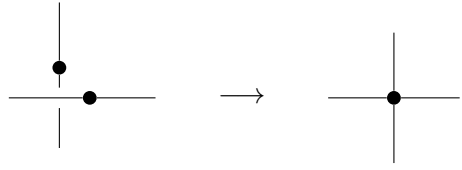
\begin{figure}[ht]
	\begin{tikzpicture}
	\node[circle,fill=black,scale=0.5] at (0-4,0) (A) {};
	\node[circle,fill=black,scale=0.5] at (0-4-0.4,0+0.4) (A') {};
	\node at (0-4-0.4,0) (A'') {};
	\node at (1-4,0) (B) {};
	\node at (0-0.4-4,1+0.4) (C) {};
	\node at (-1-0.2-4,0) (D) {};
	\node at (0-0.4-4,-1+0.2) (E) {};

	\draw (B) -- (A) -- (D);
	\draw (C) -- (A') -- (A'') -- (E);

	\node at (-2,0) {$\longrightarrow$};

	\node[circle,fill=black,scale=0.5] at (0,0) (a) {};
	\node at (1,0) (b) {};
	\node at (0,1) (c) {};
	\node at (-1,0) (d) {};
	\node at (0,-1) (e) {};
	\draw (a)--(b);
	\draw (a)--(c);
	\draw (a)--(d);
	\draw (a)--(e);
	\end{tikzpicture}
	\caption{The normalisation of a $4$-valent vertex in the toric
$1$-skeleton.}
	\label{fig:normalisation}
 \end{figure}

Let $\Delta$ denote the toric $1$-skeleton of $Y$.  Consider the partial normalisation $\Delta' \to \Delta$ constructed by normalising each $4$-valent vertex as shown in Figure~\ref{fig:normalisation}.
The partial normalisation $\Delta'$ consists of rational curves that meet at either bivalent or trivalent vertices.  Let $e$ denote an edge of $P$, and $\Gamma_e$ denote the corresponding toric curve in $X_P$.  We define two partitions of the set $\widetilde{\Gamma}_e$, as follows.   Let $p \in X_P$ be one of the endpoints of $\Gamma_e$, and $F \subset P$ the corresponding face of $P$.  Consider the part $\Delta'_p$ of the partially-normalised toric 1-skeleton $\Delta'$ that lies over $p$; this consists of the dual graph to the polyhedral subdivision (B) of $F$, partially normalised as described in Figure~\ref{fig:normalisation}.  The components of the partially-normalised dual graph define a partition of $\widetilde{\Gamma}_e$.  There is one such partition for each of the two endpoints $p$ of $\Gamma_e$: let us denote them by $\Pi_e$ and $\Pi'_e$ (the order will not matter).  The condition that Corti--Hacking--Petracci impose on their decomposition data is: for each edge $e$ of $P$ such that the dual edge $e^\star$ has lattice length $\ell(e^\star)$ equal to $1$ and for each pair $T \in \Pi$,~$T' \in \Pi'$, we have $|T \cap T'| \leq 1$.

\subsection{Computing the Betti numbers of the smoothing}

Corti--Hacking--Petracci prove:
\begin{theorem}

Let $P$ be a $3$-dimensional reflexive polytope and let $X_P$ be the toric Fano threefold associated to the spanning fan of $P$.  Fix decomposition data for $P$ as above, and let $\pi \colon Y \to X_P$ be the associated crepant toric partial resolution. Then:
\begin{enumerate}
\item $Y$ is unobstructed and smoothable;
\item If $Y_t$ is a general smoothing of $Y$, then $Y_t$ is a weak Fano threefold and the anticanonical morphism $\pi_t \colon Y_t \to X_t$,
where $X_t = \Proj R (Y_t, -K_{Y_t})$, contracts a finite number of disjoint nonsingular rational curves, each with normal bundle $\cO(-1)\oplus \cO(-1)$;
\item $X_t$ is a Fano threefold with ordinary nodes as singularities and it is a deformation of $X_P$.
\end{enumerate}
\end{theorem}
\noindent A theorem of Namikawa~\cite[Theorem~11]{Namikawa1997} now implies that $X_t$ is smoothable.  It follows that $X_P$ is smoothable.  Let $X_\eta$ denote a generic smoothing of $X_t$.  Our goal is to compute the Betti numbers of $X_\eta$.

We begin by analysing the topology of $Y_t$.  The vanishing cycles for the degeneration $Y_t \rightsquigarrow Y$ are three-dimensional spheres.  The sheaf of vanishing cycles $\phi_f$ from \eqref{short_van}, which is concentrated at the nodes of $Y$, therefore has stalk at each node equal to $\QQ$ concentrated in degree three, and the vanishing cycle exact sequence \eqref{eq:van} gives
\begin{align}
    \begin{aligned} 
    	b_0(Y_t) &= b_0(Y)  &  b_4(Y_t) &= b_2(Y)   \\
    	b_1(Y_t) &= b_1(Y)  &  b_5(Y_t) &= b_1(Y)  \\
    	b_2(Y_t) &= b_2(Y) & b_6(Y_t) &= b_0(Y) \\
    	b_3(Y_t) &= b_3(Y) - b_4(Y) + b_2(Y) + k && 
    \end{aligned} \label{eq:Y_to_Yt}
\intertext{where $k$ is the number of nodes on $Y$.  Note that $H_2(Y_t)$ is canonically identified with $H_2(Y)$; also $b_5(Y_t) = b_1(Y_t) = 0$ by Poincar\'e duality and Kodaira vanishing.  Passing from $Y_t$ to $X_\eta$ is an example of a conifold transition~\cite{CandelasGreenHubsch90}, and therefore}
    \begin{aligned} 
    	b_0(X_\eta) &= b_0(Y_t) = 1 &  b_4(X_\eta) &= b_2(Y_t) - l   \\
    	b_1(X_\eta) &= b_1(Y_t) = 0 &  b_5(X_\eta) &= b_1(Y_t) = 0  \\
    	b_2(X_\eta) &= b_2(Y_t) - l & b_6(X_\eta) &= b_0(Y_t) = 1 \\
    	b_3(X_\eta) &= b_3(Y_t) + 2m - 2l &
    \end{aligned} \label{eq:Yt_to_X}
\end{align}
where $l$ is the dimension of the subspace $L$ of $H_2(Y_t)$ spanned by the classes of curves that are contracted by $\pi_t \colon Y_t \to X_t$, and $m$ is the number of nodes on $X_t$.

In view of Theorem~\ref{thm:toric_Betti}, computing the Betti numbers of the smoothing $X_\eta$ comes down to computing the integers $k$,~$l$, and~$m$.  We have seen that $k$ is the number of quadrilaterals in the polyhedral decomposition of the boundary of $P$ determined by the decomposition data.  Corti--Hacking--Petracci conjecture the values of $l$ and $m$, as follows.  Let us identify $L \subset H_2(Y_t)$ as a subspace of $H_2(Y)$ via the canonical isomorphism $H_2(Y_t) \cong H_2(Y)$ just discussed.  Recall that an edge $e$ of $P$ determines a toric rational curve $\Gamma_e$ in $X_P$: this is the toric subvariety of $X_P$ defined by the cone over $e$.  The set $\widetilde{\Gamma}_e$ of toric curves that map dominantly to $\Gamma_e$ under the map $Y \to X_P$ is indexed by the $\ell(e)$ line segments that subdivide~$e$.  If $\ell(e) \geq 2$ then $e$ contains interior lattice points: such a lattice point $v$ then corresponds to a toric surface $S_v \subset Y$ that projects to $\Gamma_e$.  Let us denote the homology class\footnote{The fiber $F_v$ \label{fn:homology} is homologous in $Y$ to the sum of toric curves corresponding to $2$-dimensional cones in~$\Sigma$ that contain the ray spanned by $v$ and that lie entirely on one side of the hyperplane defined by the edge~$e$. The choice of side does not matter here, as the resulting sums are homologous. } of a fiber of $S_v \to \Gamma_e$ by~$F_v$. For a pair of toric curves $C$,~$C' \in \widetilde{\Gamma}_e$, let $c$,~$c'$ denote the corresponding line segments and $(c, c') \subset P$ denote the relative interior of the convex hull of $c$ and $c'$. Recall the two partitions of $\widetilde{\Gamma}_e$ defined in the discussion around Figure~\ref{fig:normalisation}, and define $n_{c,c'}$ by starting with $\ell(e)$ and subtracting one for each partition that has~$c$ and~$c'$ in the same part.  The set of exceptional curves for $\pi_t \colon Y_t \to X_t$ is conjecturally indexed \label{exceptional} by edges $e$ of $P$ such that $\ell(e) \geq 2$ and pairs of distinct elements $c$,~$c' \in \widetilde{\Gamma}_e$: it contains precisely $n_{c,c'}$ curves in the homology class
\[
\sum_{v \in (c,c')} F_v
\]
and no others.  Note that $n_{c,c'}$ here is non-negative; it can be zero.  This conjecture determines the subspace $L \subset H_2(Y_t)$ spanned by exceptional curves, and thus determines $l = \dim L$.  For $m$, suppose that the elements of the two partitions of $\widetilde{\Gamma}_e$ have sizes $a_1,a_2,\ldots$ and $b_1,b_2,\ldots$ respectively.  Set
\begin{equation} \label{eq:N_e}
N_e = \ell(e^\star) { \ell(e) \choose 2 } - \sum_i { a_i \choose 2 } - \sum_i { b_i \choose 2}.
\end{equation}
Then
\begin{equation} \label{eq:m}
m = \sum_{e \colon \ell(e) \geq 2} N_e.
\end{equation}

\subsection{Minkowski polynomials and smoothings} \label{sec:topology}

An admissible Minkowski decomposition of a three-dimensional reflexive polytope $P$ determines a Minkowski polynomial~\cite{AkhtarCoatesGalkinKasprzyk2012}.  This is a Laurent polynomial with Newton polytope $P$.  In the notation of the discussion on page~\pageref{choices}, it depends on the choices (A) of admissible Minkowski decomposition of each facet of $P$, but not on the choices (B) of regular fine mixed subdivision.

It is known that Minkowski polynomials provide mirrors for three-dimensional Fano manifolds~\cite{CoatesCortiGalkinKasprzyk2016}. As mentioned in the Introduction, if a Fano manifold $X$ is mirror to a Laurent polynomial $f$ with Newton polytope $P$, it is expected that there is a degeneration $\cX \to \Delta$ with general fiber $X$ and special fiber $X_P$.  For this expectation to be compatible with the results and conjectures of Corti--Hacking--Petracci, therefore, the Betti numbers of the Corti--Hacking--Petracci smoothings $X_\eta$ must depend only on the Minkowski polynomial.  That is, the Betti numbers must depend only on decomposition data only through the choices (A) of admissible Minkowski decompositions of facets: they must be independent of the choices (B) of regular fine mixed subdivision.  This is not obvious from the construction; we prove it in Section~\ref{sec:main_result}.

\section{Examples}
\label{sec:examples}
\subsection{Cube}

\begin{figure}[ht]
	\centering
\tdplotsetmaincoords{80}{110}
	\begin{tikzpicture}[line join=bevel,tdplot_main_coords]

\draw[gray, thick] (-1,-1,1) -- (-1,1,1) -- (1,1,1) -- (1,-1,1) -- cycle;
\draw[gray, thick] (-1,1,-1) -- (-1,1,1) -- (1,1,1) -- (1,1,-1) -- cycle;
\draw[gray, thick] (1,-1,-1) -- (1,1,-1) -- (1,1,1) -- (1,-1,1) -- cycle;

\draw[gray, thin] (-1,-1,-1) -- (-1,1,-1);
\draw[gray, thin] (-1,-1,-1) -- (1,-1,-1);
\draw[gray, thin] (-1,-1,-1) -- (-1,-1,1);

	\foreach \x in {-1,0,1} {
		\foreach \y in {-1,0,1} {
			\foreach \z in {-1,0,1} {
				\pgfmathsetmacro\mycolour{ifthenelse(\x*\y*\z==0,"gray","black")};
				\filldraw[\mycolour] (\x,\y,\z) circle (1pt);
			}
		}
	}
\end{tikzpicture}
\hspace{2.4cm}
\begin{tikzpicture}[line join=bevel,tdplot_main_coords,scale=1.2]

\draw[black, semithick] (0,0,1) -- (1,0,0) -- (0,1,0) -- cycle;
\draw[black, semithick] (0,0,1) -- (0,-1,0) -- (1,0,0) -- cycle;

\draw[black, semithick] (0,0,-1) -- (1,0,0) -- (0,1,0) -- cycle;
\draw[black, semithick] (0,0,-1) -- (0,-1,0) -- (1,0,0) -- cycle;

\draw[black, ultra thin] (-1,0,0) -- (0,0,1);
\draw[black, ultra thin] (-1,0,0) -- (0,0,-1);
\draw[black, ultra thin] (-1,0,0) -- (0,1,0);
\draw[black, ultra thin] (-1,0,0) -- (0,-1,0);

\end{tikzpicture}

\caption{The cube $P$ and a schematic picture of the toric variety $X_P$}
\label{fig:cube_and_octahedron}
\end{figure}
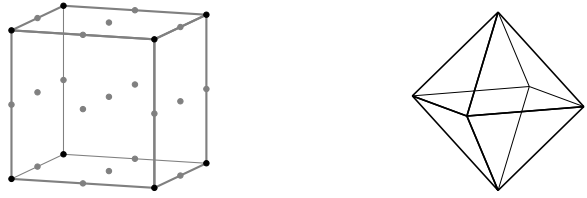

Consider the cube $P$ centered at the origin with vertices $(\pm 1,\pm 1, \pm 1)$.  This has six non-simplicial facets and twelve edges of length two; thus the toric variety $X_P$ defined by the spanning fan of $P$ has six singular points and twelve curves of transverse~$A_1$ singularities. These are arranged as on the right-hand side of Figure~\ref{fig:cube_and_octahedron}, with the singular points at the vertices of the octahedron and the singular curves as the edges.

\begin{figure}[ht]
\centering
\begin{tikzpicture}

\draw[gray, thick] (-1,-1) -- (-1,1);
\draw[gray, thick] (-1,-1) -- (1,-1);
\draw[gray, thick] (1,1) -- (-1,1);
\draw[gray, thick] (1,1) -- (1,-1);
\draw[gray, thick] (-1,0) -- (1,0);
\draw[gray, thick] (0,-1) -- (0,1);

\filldraw[black] (-1,-1) circle (1pt);
\filldraw[black] (-1,0) circle (1pt);
\filldraw[black] (-1,1) circle (1pt);
\filldraw[black] (0,-1) circle (1pt);
\filldraw[black] (0,0) circle (1pt);
\filldraw[black] (0,1) circle (1pt);
\filldraw[black] (1,-1) circle (1pt);
\filldraw[black] (1,0) circle (1pt);
\filldraw[black] (1,1) circle (1pt);

\end{tikzpicture}
\caption{The unique fine mixed subdivision of the facet $F$}
\label{fig:cube_facet}
\end{figure}
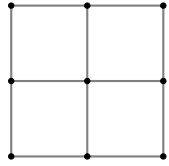

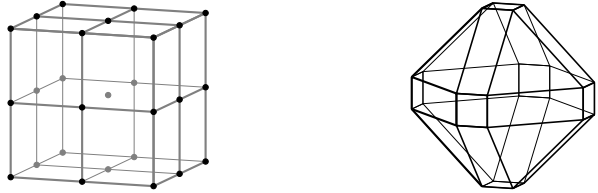
\begin{figure}[ht]
\tdplotsetmaincoords{80}{110}
	\begin{tikzpicture}[line join=bevel,tdplot_main_coords]

\draw[gray, thick] (-1,-1,1) -- (-1,1,1) -- (1,1,1) -- (1,-1,1) -- cycle;
\draw[gray, thick] (-1,1,-1) -- (-1,1,1) -- (1,1,1) -- (1,1,-1) -- cycle;
\draw[gray, thick] (1,-1,-1) -- (1,1,-1) -- (1,1,1) -- (1,-1,1) -- cycle;

\draw[gray, ultra thin] (-1,-1,-1) -- (-1,1,-1);
\draw[gray, ultra thin] (-1,-1,-1) -- (1,-1,-1);
\draw[gray, ultra thin] (-1,-1,-1) -- (-1,-1,1);

\draw[gray, thick] (0,-1,1) -- (0,1,1);
\draw[gray, thick] (-1,0,1) -- (1,0,1);
\draw[gray, ultra thin] (0,-1,-1) -- (0,1,-1);
\draw[gray, ultra thin] (-1,0,-1) -- (1,0,-1);

\draw[gray, thick] (0,1,-1) -- (0,1,1);
\draw[gray, thick] (-1,1,0) -- (1,1,0);
\draw[gray, ultra thin] (0,-1,-1) -- (0,-1,1);
\draw[gray, ultra thin] (-1,-1,0) -- (1,-1,0);

\draw[gray, thick] (1,-1,0) -- (1,1,0);
\draw[gray, thick] (1,0,-1) -- (1,0,1);
\draw[gray, ultra thin] (-1,-1,0) -- (-1,1,0);
\draw[gray, ultra thin] (-1,0,-1) -- (-1,0,1);

	\foreach \x in {-1,0,1} {
		\foreach \y in {-1,0,1} {
			\foreach \z in {-1,0,1} {
				\pgfmathsetmacro\mycolour{ifthenelse(\x==1||\y==1||\z==1,"black","gray")};
				\filldraw[\mycolour] (\x,\y,\z) circle (1pt);
			}
		}
	}

\end{tikzpicture}
\hspace{2.4cm}
\begin{tikzpicture}[line join=bevel,tdplot_main_coords,scale=1.2]

\draw[black, semithick] (1,0.18,0.18) -- (1,0.18,-0.18) -- (1,-0.18,-0.18) -- (1,-0.18,0.18) -- cycle;
\draw[black, semithick] (0.18,1,0.18) -- (0.18,1,-0.18) -- (-0.18,1,-0.18) -- (-0.18,1,0.18) -- cycle;
\draw[black, semithick] (0.18,0.18,1) -- (0.18,-0.18,1) -- (-0.18,-0.18,1) -- (-0.18,0.18,1) -- cycle;
\draw[black, ultra thin] (-1,0.18,0.18) -- (-1,0.18,-0.18) -- (-1,-0.18,-0.18) -- (-1,-0.18,0.18) -- cycle;
\draw[black, ultra thin] (0.18,-1,0.18) -- (0.18,-1,-0.18) -- (-0.18,-1,-0.18) -- (-0.18,-1,0.18) -- cycle;
\draw[black, ultra thin] (0.18,0.18,-1) -- (0.18,-0.18,-1) -- (-0.18,-0.18,-1) -- (-0.18,0.18,-1) -- cycle;

\draw[black, semithick] (1,0.18,0.18) -- (0.18,0.18,1) -- (0.18,-0.18,1) -- (1,-0.18,0.18) -- cycle;
\draw[black, semithick] (1,0.18,0.18) -- (0.18,1,0.18) -- (0.18,1,-0.18) -- (1,0.18,-0.18) -- cycle;
\draw[black, semithick] (0.18,1,0.18) -- (0.18,0.18,1) -- (-0.18,0.18,1) -- (-0.18,1,0.18) -- cycle;
\draw[black, semithick] (1,0.18,-0.18) -- (0.18,0.18,-1) -- (0.18,-0.18,-1) -- (1,-0.18,-0.18) -- cycle;

\draw[black, semithick] (0.18,0.18,-1) -- (-0.18,0.18,-1) -- (-0.18,1,-0.18) -- (0.18,1,-0.18) -- cycle;

\draw[black, ultra thin] (-0.18,1,0.18) -- (-1,0.18,0.18) -- (-1,0.18,-0.18) -- (-0.18,1,-0.18) -- cycle;
\draw[black, ultra thin] (-1,0.18,0.18) -- (-0.18,0.18,1) -- (-0.18,-0.18,1) -- (-1,-0.18,0.18) -- cycle;
\draw[black, ultra thin] (-1,-0.18,0.18) -- (-0.18,-1,0.18) -- (-0.18,-1,-0.18) -- (-1,-0.18,-0.18) -- cycle;

\draw[black, ultra thin] (-0.18,-1,0.18) -- (0.18,-1,0.18) -- (0.18,-0.18,1) -- (-0.18,-0.18,1) -- cycle;
\draw[black, thick] (0.18,-1,0.18) -- (0.18,-0.18,1) -- (-0.18,-0.18,1);

\draw[black, ultra thin] (-0.18,-1,-0.18) -- (0.18,-1,-0.18) -- (0.18,-0.18,-1) -- (-0.18,-0.18,-1) -- cycle;
\draw[black, thick] (0.18,-1,-0.18) -- (0.18,-0.18,-1);

\draw[black, thick] (1,-0.18,-0.18) --  (1,-0.18,0.18) -- (0.18,-1,0.18) -- (0.18,-1,-0.18) -- cycle;

\draw[black, ultra thin] (-0.18,-0.18,-1) -- (-1,-0.18,-0.18) -- (-1,0.18,-0.18) -- (-0.18,0.18,-1) -- cycle;
\end{tikzpicture}

\caption{The polyhedral subdivision and a schematic picture of $Y$}
\label{fig:cube_decomposition}
\end{figure}

 Each facet $F$ of $P$ is a square with side-length two; this has a unique admissible Minkowski decomposition, as a Minkowski sum of four line segments, which in turn leads to the unique fine mixed subdivision of $F$ shown in Figure~\ref{fig:cube_facet}.
The fan $\Sigma$ that defines the partial resolution $Y$ of $X_P$ is therefore obtained by taking cones over the polyhedral decomposition of the boundary of $P$ shown in Figure~\ref{fig:cube_decomposition}.
The variety $Y$ contains 24 ordinary double points and 48 toric curves, arranged as on the right-hand side of Figure~\ref{fig:cube_decomposition}, with the singular points at the vertices and the toric curves (along which $Y$ is non-singular) as the edges.  Furthermore each edge in the dual polygon $P^\star$ has length~$1$, and each partition $\Pi_e$ of $\widetilde{\Gamma}_e$ is into singleton sets, so the polyhedral decomposition satisfies the conditions to be decomposition data.

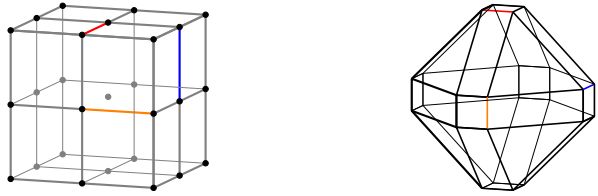
\begin{figure}[ht]
\tdplotsetmaincoords{80}{110}
	\begin{tikzpicture}[line join=bevel,tdplot_main_coords]

\draw[gray, thick] (-1,-1,1) -- (-1,1,1) -- (1,1,1) -- (1,-1,1) -- cycle;
\draw[gray, thick] (-1,1,-1) -- (-1,1,1) -- (1,1,1) -- (1,1,-1) -- cycle;
\draw[gray, thick] (1,-1,-1) -- (1,1,-1) -- (1,1,1) -- (1,-1,1) -- cycle;

\draw[gray, ultra thin] (-1,-1,-1) -- (-1,1,-1);
\draw[gray, ultra thin] (-1,-1,-1) -- (1,-1,-1);
\draw[gray, ultra thin] (-1,-1,-1) -- (-1,-1,1);

\draw[gray, thick] (0,-1,1) -- (0,1,1);
\draw[gray, thick] (-1,0,1) -- (0,0,1);
\draw[red, thick] (0,0,1) -- (1,0,1);
\draw[gray, ultra thin] (0,-1,-1) -- (0,1,-1);
\draw[gray, ultra thin] (-1,0,-1) -- (1,0,-1);

\draw[gray, thick] (0,1,-1) -- (0,1,0);
\draw[blue, thick] (0,1,0) -- (0,1,1);
\draw[gray, thick] (-1,1,0) -- (1,1,0);
\draw[gray, ultra thin] (0,-1,-1) -- (0,-1,1);
\draw[gray, ultra thin] (-1,-1,0) -- (1,-1,0);

\draw[gray, thick] (1,-1,0) -- (1,0,0);
\draw[orange, thick] (1,0,0) -- (1,1,0);
\draw[gray, thick] (1,0,-1) -- (1,0,1);
\draw[gray, ultra thin] (-1,-1,0) -- (-1,1,0);
\draw[gray, ultra thin] (-1,0,-1) -- (-1,0,1);

	\foreach \x in {-1,0,1} {
		\foreach \y in {-1,0,1} {
			\foreach \z in {-1,0,1} {
				\pgfmathsetmacro\mycolour{ifthenelse(\x==1||\y==1||\z==1,"black","gray")};
				\filldraw[\mycolour] (\x,\y,\z) circle (1pt);
			}
		}
	}

\end{tikzpicture}
\hspace{2.4cm}
\begin{tikzpicture}[line join=bevel,tdplot_main_coords,scale=1.2]

\draw[black, semithick] (1,0.18,-0.18) -- (1,-0.18,-0.18) -- (1,-0.18,0.18) -- (1,0.18,0.18);
\draw[black, semithick] (0.18,1,0.18) -- (0.18,1,-0.18) -- (-0.18,1,-0.18) -- (-0.18,1,0.18);
\draw[black, semithick] (0.18,-0.18,1) -- (-0.18,-0.18,1) -- (-0.18,0.18,1);
\draw[black, ultra thin] (-1,0.18,0.18) -- (-1,0.18,-0.18) -- (-1,-0.18,-0.18) -- (-1,-0.18,0.18) -- cycle;
\draw[black, ultra thin] (0.18,-1,0.18) -- (0.18,-1,-0.18) -- (-0.18,-1,-0.18) -- (-0.18,-1,0.18) -- cycle;
\draw[black, ultra thin] (0.18,0.18,-1) -- (0.18,-0.18,-1) -- (-0.18,-0.18,-1) -- (-0.18,0.18,-1) -- cycle;

\draw[black, semithick] (0.18,-0.18,1) -- (1,-0.18,0.18) -- (1,0.18,0.18) -- (0.18,0.18,1);
\draw[red, semithick] (0.18,0.18,1) -- (0.18,-0.18,1);
\draw[black, semithick] (1,0.18,0.18) -- (0.18,1,0.18) -- (0.18,1,-0.18) -- (1,0.18,-0.18);
\draw[orange, semithick] (1,0.18,0.18) -- (1,0.18,-0.18);
\draw[black, semithick] (0.18,1,0.18) -- (0.18,0.18,1) -- (-0.18,0.18,1) -- (-0.18,1,0.18);
\draw[blue, semithick] (0.18,1,0.18) -- (-0.18,1,0.18);
\draw[black, semithick] (1,0.18,-0.18) -- (0.18,0.18,-1) -- (0.18,-0.18,-1) -- (1,-0.18,-0.18) -- cycle;
\draw[black, semithick] (0.18,0.18,-1) -- (-0.18,0.18,-1) -- (-0.18,1,-0.18) -- (0.18,1,-0.18) -- cycle;

\draw[black, ultra thin] (-0.18,1,0.18) -- (-1,0.18,0.18) -- (-1,0.18,-0.18) -- (-0.18,1,-0.18) -- cycle;
\draw[black, ultra thin] (-1,0.18,0.18) -- (-0.18,0.18,1) -- (-0.18,-0.18,1) -- (-1,-0.18,0.18) -- cycle;
\draw[black, ultra thin] (-1,-0.18,0.18) -- (-0.18,-1,0.18) -- (-0.18,-1,-0.18) -- (-1,-0.18,-0.18) -- cycle;

\draw[black, ultra thin] (-0.18,-1,0.18) -- (0.18,-1,0.18) -- (0.18,-0.18,1) -- (-0.18,-0.18,1) -- cycle;
\draw[black, thick] (0.18,-1,0.18) -- (0.18,-0.18,1) -- (-0.18,-0.18,1);

\draw[black, ultra thin] (-0.18,-1,-0.18) -- (0.18,-1,-0.18) -- (0.18,-0.18,-1) -- (-0.18,-0.18,-1) -- cycle;
\draw[black, thick] (0.18,-1,-0.18) -- (0.18,-0.18,-1);

\draw[black, thick] (1,-0.18,-0.18) --  (1,-0.18,0.18) -- (0.18,-1,0.18) -- (0.18,-1,-0.18) -- cycle;

\draw[black, ultra thin] (-0.18,-0.18,-1) -- (-1,-0.18,-0.18) -- (-1,0.18,-0.18) -- (-0.18,0.18,-1) -- cycle;
\end{tikzpicture}

\caption{Toric curves that generate the subspace of exceptional curves}
\label{fig:cube_generators}
\end{figure}

Applying Theorem~\ref{thm:toric_Betti} gives
\begin{align*}
	b_0(Y) &= 1  &  b_4(Y) &= 23   \\
	b_1(Y) &= 0  &  b_5(Y) &= 0  \\
	b_2(Y) &= 4 & b_6(Y) &= 1 \\
	b_3(Y) &= 5 &
\end{align*}
There are $k=24$ quadrilaterals in the polyhedral subdivision of the boundary of~$P$, and the discussion in Section~\ref{sec:topology} yields
\begin{align*}
	b_0(Y_t) &= 1  &  b_4(Y_t) &= 4   \\
	b_1(Y_t) &= 0  &  b_5(Y_t) &= 0  \\
	b_2(Y_t) &= 4 & b_6(Y_t) &= 1. \\
	b_3(Y_t) &= 10 &
\end{align*}
From the description of the homology classes $F_v$ in footnote~\ref{fn:homology} and~\cite[Proposition~2.1]{FultonSturmfels1997}, we see that $l=3$, with generators for the subspace $L \subset H_2(Y_t) \cong H_2(Y)$ of classes of exceptional curves as shown in Figure~\ref{fig:cube_generators}.

\noindent  Equation \eqref{eq:m} gives $m=12$, and we find:
\begin{align*}
	b_0(X_\eta) &= 1  &  b_4(X_\eta) &= 1   \\
	b_1(X_\eta) &= 0  &  b_5(X_\eta) &= 0  \\
	b_2(X_\eta) &= 1 & b_6(X_\eta) &= 1 \\
	b_3(X_\eta) &= 28 &
\end{align*}
The only Fano 3-fold with these Betti numbers is $V_8$.  Thus $X_\eta$ is isomorphic to $V_8$, which is consistent with the fact that the Minkowski polynomial
\[
f = \frac{(1+x)^2 (1+y)^2 (1+z)^2}{xyz} - 8
\]
defined by our decomposition data for $P$ is a mirror to $V_8$.

\subsection{A singular toric variety with two different smoothings}

\begin{figure}[ht]
\tdplotsetmaincoords{70}{106}
\begin{tikzpicture}[line join=bevel,tdplot_main_coords]

	\draw[gray, thick] (0,0,1) -- (0,1,-1) -- (1,1,-1) -- cycle;
	\draw[gray, thick] (0,0,1) -- (1,1,-1) -- (1,0,-1) -- cycle;
	\draw[gray, thick] (0,0,1) -- (1,0,-1) -- (0,-1,-1) -- cycle;
	\draw[gray, thin] (0,0,1) -- (0,-1,-1) -- (-1,-1,-1) -- cycle;
	\draw[gray, thin] (0,0,1) -- (-1,-1,-1) -- (-1,0,-1) -- cycle;
	\draw[gray, thin] (0,0,1) -- (-1,0,-1) -- (0,1,-1) -- cycle;

	\filldraw[gray] (0,0,0) circle (1pt);
	\filldraw[gray] (0,0,-1) circle (1pt);
	\filldraw[black] (0,0,1) circle (1pt);
	\filldraw[black] (0,1,-1) circle (1pt);
	\filldraw[black] (1,1,-1) circle (1pt);
	\filldraw[black] (1,0,-1) circle (1pt);
	\filldraw[black] (0,-1,-1) circle (1pt);
	\filldraw[gray] (-1,-1,-1) circle (1pt);
	\filldraw[gray] (-1,0,-1) circle (1pt);
\end{tikzpicture}
\hspace{2.4cm}
\begin{tikzpicture}[line join=bevel,tdplot_main_coords,scale=1.2]

\draw[black, ultra thin] (0,0,1) -- (-2,0,-1) -- (0,-2,-1) -- cycle;
\draw[black, semithick] (0,0,1) -- (0,-2,-1) -- (2,-2,-1) -- cycle;
\draw[black, semithick] (0,0,1) -- (2,-2,-1) -- (2,0,-1) -- cycle;
\draw[black, semithick] (0,0,1) -- (2,0,-1) -- (0,2,-1) -- cycle;
\draw[black, semithick] (0,0,1) -- (0,2,-1) -- (-2,2,-1) -- cycle;
\draw[black, ultra thin] (0,0,1) -- (-2,2,-1) -- (-2,0,-1) -- cycle;

\end{tikzpicture}
\caption{The polytope $P$ and a schematic picture of the toric variety $X_P$}
\label{fig:122}
\end{figure}
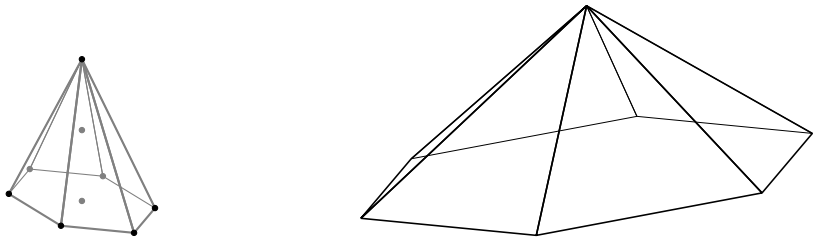

Consider the three-dimensional polytope, pictured in Figure~\ref{fig:122}, with vertices
\[
\text{$(0,0,1)$, $(0,1,-1)$, $(1,1,-1)$, $(1,0,-1)$, $(0,-1,-1)$, $(-1,-1,-1)$, $(-1,0,-1)$.}
\]
This polytope $P$ is reflexive. It has six facets that are standard simplices, one non-simplicial facet (a hexagon), and 12 edges of length~$1$.  Thus the toric variety $X_P$ defined by the spanning fan of $P$ contains six non-singular toric points, a unique singular point (at which the singularity is a cone over the del~Pezzo surface of degree~$6$), and 12 non-singular toric curves.  These are arranged as on the right-hand side of Figure~\ref{fig:122}, with the non-singular toric points at the $3$-valent vertices, the singular point at the $6$-valent vertex, and the toric curves as the edges.

\begin{figure}[ht]
\centering
\begin{tikzpicture}
\draw[gray, thick] (0,0) -- (0,-1) -- (1,0) -- (1,1) -- cycle;
\draw[gray, thick] (0,0) -- (0,-1) -- (-1,-1) -- (-1,0) -- cycle;
\draw[gray, thick] (0,0) -- (-1,0) -- (0,1) -- (1,1) -- cycle;
\filldraw[black] (-1,-1) circle (1pt);
\filldraw[black] (-1,0) circle (1pt);
\filldraw[black] (0,-1) circle (1pt);
\filldraw[black] (0,0) circle (1pt);
\filldraw[black] (0,1) circle (1pt);
\filldraw[black] (1,0) circle (1pt);
\filldraw[black] (1,1) circle (1pt);
\end{tikzpicture}
\hspace{3cm}
\begin{tikzpicture}
\draw[gray, thick] (0,0) -- (-1,-1) -- (-1,0) -- cycle;
\draw[gray, thick] (0,0) -- (1,0) -- (1,1) -- cycle;
\draw[gray, thick] (0,0) -- (-1,0) -- (0,1) -- (1,1) -- cycle;
\draw[gray, thick] (0,0) -- (1,0) -- (0,-1) -- (-1,-1) -- cycle;
\filldraw[black] (-1,-1) circle (1pt);
\filldraw[black] (-1,0) circle (1pt);
\filldraw[black] (0,-1) circle (1pt);
\filldraw[black] (0,0) circle (1pt);
\filldraw[black] (0,1) circle (1pt);
\filldraw[black] (1,0) circle (1pt);
\filldraw[black] (1,1) circle (1pt);
\end{tikzpicture}
\caption{Two fine mixed subdivisions of the facet $F$}
\label{fig:122_facet}
\end{figure}
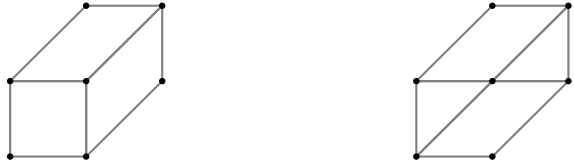

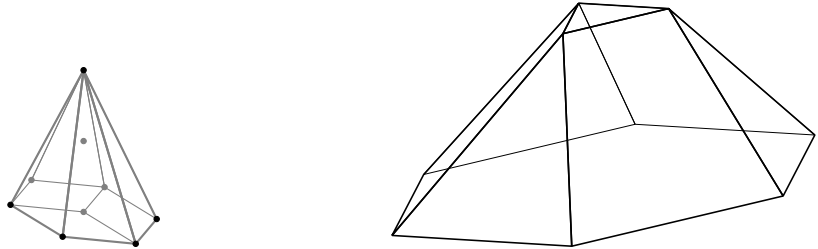
\begin{figure}[ht]
\tdplotsetmaincoords{70}{106}
\begin{tikzpicture}[line join=bevel,tdplot_main_coords]

	\draw[gray, thick] (0,0,1) -- (0,1,-1) -- (1,1,-1) -- cycle;
	\draw[gray, thick] (0,0,1) -- (1,1,-1) -- (1,0,-1) -- cycle;
	\draw[gray, thick] (0,0,1) -- (1,0,-1) -- (0,-1,-1) -- cycle;
	\draw[gray, thin] (0,0,1) -- (0,-1,-1) -- (-1,-1,-1) -- cycle;
	\draw[gray, thin] (0,0,1) -- (-1,-1,-1) -- (-1,0,-1) -- cycle;
	\draw[gray, thin] (0,0,1) -- (-1,0,-1) -- (0,1,-1) -- cycle;

	\draw[gray, thin] (0,0,-1) -- (0,-1,-1);
	\draw[gray, thin] (0,0,-1) -- (1,1,-1);
	\draw[gray, thin] (0,0,-1) -- (-1,0,-1);

	\filldraw[gray] (0,0,0) circle (1pt);
	\filldraw[gray] (0,0,-1) circle (1pt);
	\filldraw[black] (0,0,1) circle (1pt);
	\filldraw[black] (0,1,-1) circle (1pt);
	\filldraw[black] (1,1,-1) circle (1pt);
	\filldraw[black] (1,0,-1) circle (1pt);
	\filldraw[black] (0,-1,-1) circle (1pt);
	\filldraw[gray] (-1,-1,-1) circle (1pt);
	\filldraw[gray] (-1,0,-1) circle (1pt);
\end{tikzpicture}
\hspace{2.4cm}
\tdplotsetmaincoords{70}{100}
\begin{tikzpicture}[line join=bevel,tdplot_main_coords,scale=1.2]

\draw[black, semithick] (-1/3,-1/3,1) -- (2/3,-1/3,1) -- (-1/3,2/3,1) -- cycle;
\draw[black, ultra thin] (-1/3,-1/3,1) -- (-2,0,-1) -- (0,-2,-1) -- cycle;
\draw[black, semithick] (2/3,-1/3,1) -- (-1/3,-1/3,1) -- (0,-2,-1) -- (2,-2,-1) -- cycle;
\draw[black, semithick] (2/3,-1/3,1) -- (2,-2,-1) -- (2,0,-1) -- cycle;
\draw[black, semithick] (-1/3,2/3,1) -- (2/3,-1/3,1) -- (2,0,-1) -- (0,2,-1) -- cycle;
\draw[black, semithick] (-1/3,2/3,1) -- (0,2,-1) -- (-2,2,-1) -- cycle;
\draw[black, ultra thin] (-1/3,-1/3,1) -- (-1/3,2/3,1) -- (-2,2,-1) -- (-2,0,-1) -- cycle;

\end{tikzpicture}
\caption{The first polyhedral subdivision and a schematic picture of $Y$}
\label{fig:122_first_subdivision}
\end{figure}

The hexagonal facet $F$ of $P$ admits two Minkowski decompositions: as the sum of three line segments, and as the sum of two triangles.  Up to automorphism, these each give rise to a unique fine mixed subdivision of $F$, as shown in Figure~\ref{fig:122_facet}.
Consider first the left-most fine mixed subdivision in Figure~\ref{fig:122_facet}.  This leads to the polyhedral subdivision of the boundary of $P$ shown in Figure~\ref{fig:122_first_subdivision}.
The fan $\Sigma$ given by taking cones over this polyhedral subdivision defines a toric partial resolution $Y$ of $X_P$.  This variety $Y$ has three ordinary double points, six non-singular toric points, and 15 non-singular toric curves, arranged as on the right-hand side of Figure~\ref{fig:122_first_subdivision}: the ordinary double points are the $4$-valent vertices, the non-singular points are the $3$-valent vertices, and the toric curves are the edges.  The polyhedral decomposition satisfies the conditions to be decomposition data.

Applying Theorem~\ref{thm:toric_Betti} gives
\begin{align*}
	b_0(Y) &= 1  &  b_4(Y) &= 5  \\
	b_1(Y) &= 0  &  b_5(Y) &= 0  \\
	b_2(Y) &= 2 & b_6(Y) &= 1 \\
	b_3(Y) &= 0 &
\end{align*}
and since there are $k=3$ quadrilaterals in the polyhedral subdivision of the boundary of $P$, we have
\begin{align*}
	b_0(Y_t) &= 1  &  b_4(Y_t) &= 2   \\
	b_1(Y_t) &= 0  &  b_5(Y_t) &= 0  \\
	b_2(Y_t) &= 2 & b_6(Y_t) &= 1. \\
	b_3(Y_t) &= 0 &
\end{align*}
The discussion on page~\pageref{exceptional} implies that, assuming the conjectures of Corti--Hacking--Petracci, there are no exceptional curves in $Y_t$ and so the smoothings $X_\eta$ of $X_P$ and $Y_t$ of $Y$ are isomorphic.  Thus
\begin{align*}
	b_0(X_\eta) &= 1  &  b_4(X_\eta) &= 2   \\
	b_1(X_\eta) &= 0  &  b_5(X_\eta) &= 0  \\
	b_2(X_\eta) &= 2 & b_6(X_\eta) &= 1. \\
	b_3(X_\eta) &= 0 &
\end{align*}
These are the Betti numbers of the hypersurface $W_{1,1}$ of bidegree~$(1,1)$ in $\PP^2 \times \PP^2$, which is consistent with the fact that the Minkowski polynomial
\[
f = \frac{(1+x)(1+y)(1+xy)}{xyz} + z
\]
defined by our decomposition data for $P$ is a mirror to $W_{1,1}$.

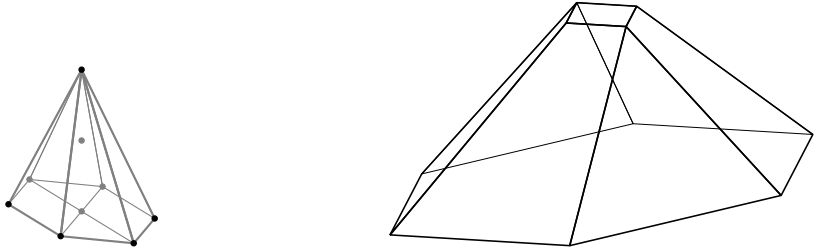
\begin{figure}[ht]
\tdplotsetmaincoords{70}{106}
\begin{tikzpicture}[line join=bevel,tdplot_main_coords]

	\draw[gray, thick] (0,0,1) -- (0,1,-1) -- (1,1,-1) -- cycle;
	\draw[gray, thick] (0,0,1) -- (1,1,-1) -- (1,0,-1) -- cycle;
	\draw[gray, thick] (0,0,1) -- (1,0,-1) -- (0,-1,-1) -- cycle;
	\draw[gray, thin] (0,0,1) -- (0,-1,-1) -- (-1,-1,-1) -- cycle;
	\draw[gray, thin] (0,0,1) -- (-1,-1,-1) -- (-1,0,-1) -- cycle;
	\draw[gray, thin] (0,0,1) -- (-1,0,-1) -- (0,1,-1) -- cycle;

	\draw[gray, thin] (0,0,-1) -- (1,0,-1);
	\draw[gray, thin] (0,0,-1) -- (1,1,-1);
	\draw[gray, thin] (0,0,-1) -- (-1,0,-1);
	\draw[gray, thin] (0,0,-1) -- (-1,-1,-1);

	\filldraw[gray] (0,0,0) circle (1pt);
	\filldraw[gray] (0,0,-1) circle (1pt);
	\filldraw[black] (0,0,1) circle (1pt);
	\filldraw[black] (0,1,-1) circle (1pt);
	\filldraw[black] (1,1,-1) circle (1pt);
	\filldraw[black] (1,0,-1) circle (1pt);
	\filldraw[black] (0,-1,-1) circle (1pt);
	\filldraw[gray] (-1,-1,-1) circle (1pt);
	\filldraw[gray] (-1,0,-1) circle (1pt);
\end{tikzpicture}
\hspace{2.4cm}
\tdplotsetmaincoords{70}{100}
\begin{tikzpicture}[line join=bevel,tdplot_main_coords,scale=1.2]

\draw[black, semithick] (-1/3,-1/3,1) -- (-1/3,1/3,1) -- (1/3,1/3,1) -- (1/3,-1/3,1) -- cycle;
\draw[black, ultra thin] (-1/3,-1/3,1) -- (-2,0,-1) -- (0,-2,-1) -- cycle;
\draw[black, semithick] (1/3,-1/3,1) -- (-1/3,-1/3,1) -- (0,-2,-1) -- (2,-2,-1) -- cycle;
\draw[black, semithick] (1/3,1/3,1) -- (1/3,-1/3,1) -- (2,-2,-1) -- (2,0,-1) -- cycle;
\draw[black, semithick] (1/3,1/3,1) -- (2,0,-1) -- (0,2,-1) -- cycle;
\draw[black, semithick] (-1/3,1/3,1) -- (1/3,1/3,1) -- (0,2,-1) -- (-2,2,-1) -- cycle;
\draw[black, ultra thin] (-1/3,-1/3,1) -- (-1/3,1/3,1) -- (-2,2,-1) -- (-2,0,-1) -- cycle;

\end{tikzpicture}
\caption{The second polyhedral subdivision and a schematic picture of $Y$}
\label{fig:122_second_subdivision}
\end{figure}

One could instead consider the right-most fine mixed subdivision in Figure~\ref{fig:122_facet}.  This leads to the polyhedral subdivision of the boundary of $P$ shown in Figure~\ref{fig:122_second_subdivision}; again this satisfies the conditions to be decomposition data.
This time the toric partial resolution $Y$ has two ordinary double points, eight non-singular toric points, and 16 non-singular toric curves.  These are arranged as on the right-hand side of Figure~\ref{fig:122_first_subdivision}, with the ordinary double points as the $4$-valent vertices, the non-singular points as the $3$-valent vertices, and the toric curves as the edges.

Theorem~\ref{thm:toric_Betti} yields
\begin{align*}
	b_0(Y) &= 1  &  b_4(Y) &= 5  \\
	b_1(Y) &= 0  &  b_5(Y) &= 0  \\
	b_2(Y) &= 3 & b_6(Y) &= 1 \\
	b_3(Y) &= 0 &
\end{align*}
and since there are $k=2$ quadrilaterals in the polyhedral subdivision of the boundary of $P$, we have
\begin{align*}
	b_0(Y_t) &= 1  &  b_4(Y_t) &= 3   \\
	b_1(Y_t) &= 0  &  b_5(Y_t) &= 0  \\
	b_2(Y_t) &= 3 & b_6(Y_t) &= 1. \\
	b_3(Y_t) &= 0 &
\end{align*}
Once again, there are (conjecturally) no exceptional curves in $Y_t$ and so the smoothings $X_\eta$ of $X_P$ and $Y_t$ of $Y$ are isomorphic.  Thus
\begin{align*}
	b_0(X_\eta) &= 1  &  b_4(X_\eta) &= 3   \\
	b_1(X_\eta) &= 0  &  b_5(X_\eta) &= 0  \\
	b_2(X_\eta) &= 3 & b_6(X_\eta) &= 1 \\
	b_3(X_\eta) &= 0 &
\end{align*}
These are the Betti numbers of $\PP^1 \times \PP^1 \times \PP^1$, which is consistent with the fact that the Minkowski polynomial
\[
f = \frac{(1+x + xy)(1+y+xy)}{xyz} + z
\]
defined by our decomposition data for $P$ is a mirror to $\PP^1 \times \PP^1 \times \PP^1$.

\subsection{An example with transverse $A_2$ singularities}
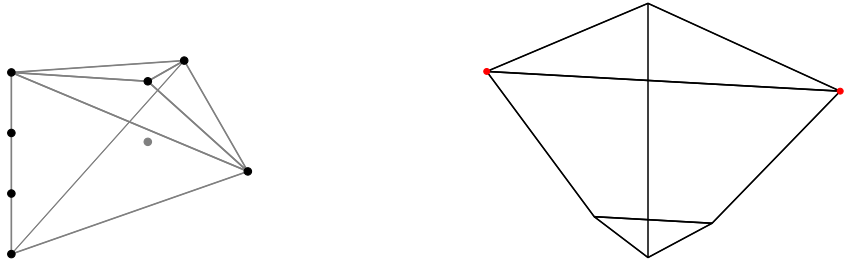
\begin{figure}[ht]
\tdplotsetmaincoords{35}{20}
\begin{tikzpicture}[line join=bevel,tdplot_main_coords,scale=1.4]

	\draw[gray, ultra thin] (1,0,0) -- (0,1,0) -- (-1,-1,-1) -- cycle;
	\draw[gray, semithick] (1,0,0) -- (0,1,0) -- (0,0,1) -- cycle;
	\draw[gray, semithick] (1,0,0) -- (0,0,1) -- (-1,-1,2) -- cycle;
	\draw[gray, ultra thin] (0,1,0) -- (-1,-1,-1) -- (-1,-1,2) -- cycle;
	\draw[gray, semithick] (0,1,0) -- (0,0,1) -- (-1,-1,2) -- cycle;
	\draw[gray, semithick] (1,0,0) -- (-1,-1,-1) -- (-1,-1,2) -- cycle;

	\filldraw[black] (1,0,0) circle (1pt);
	\filldraw[black] (0,1,0) circle (1pt);
	\filldraw[black] (0,0,1) circle (1pt);
	\filldraw[black] (-1,-1,-1) circle (1pt);
	\filldraw[black] (-1,-1,0) circle (1pt);
	\filldraw[black] (-1,-1,1) circle (1pt);
	\filldraw[black] (-1,-1,2) circle (1pt);
	\filldraw[gray] (0,0,0) circle (1pt);
\end{tikzpicture}
\hspace{2.4cm}
\tdplotsetmaincoords{50}{320}
\begin{tikzpicture}[line join=bevel,tdplot_main_coords,scale=1.1]

	\draw[black, ultra thin] (-1,-1,-1) -- (-1,-1,3) -- (2,-1,0) -- (0,-1,-1) -- cycle;
	\draw[black, ultra thin] (-1,-1,3) -- (-1,-1,-1) -- (-1,0,-1) -- (-1,2,0) -- cycle;
	\draw[black, semithick] (-1,0,-1) -- (0,-1,-1) -- (2,-1,0) -- (-1,2,0) -- cycle;
	\draw[black, semithick] (-1,-1,-1) -- (-1,0,-1) -- (0,-1,-1) -- cycle;
	\draw[black, semithick] (-1,-1,3) -- (2,-1,0) -- (-1,2,0) -- cycle;

	\filldraw[red] (2,-1,0) circle (1pt);
	\filldraw[red] (-1,2,0) circle (1pt);

\end{tikzpicture}
\caption{The polytope $P$ and a schematic picture of the toric variety $X_P$}
\label{fig:34}
\end{figure}

Consider the three-dimensional reflexive polytope $P$, pictured in Figure~\ref{fig:34}, with vertices
\[
\text{$(1,0,0)$, $(0,1,0)$, $(0,0,1)$, $(-1,-1,-1)$, $(-1,-1,2)$.}
\]
This has four facets that are standard simplices, two non-standard simplicial facets, eight edges of length~$1$, and one edge of length 3.  Thus the toric variety $X_P$ defined by the spanning fan of $P$ contains four non-singular toric points, two orbifold points, eight non-singular toric curves, and one toric curve with transverse~$A_2$ singularities.  These are arranged as on the right-hand side of Figure~\ref{fig:34}, with the toric points at the vertices, the orbifold points indicated in red, and the toric curves as the edges.  The curve of $A_2$ singularities is the edge between the two orbifold points.

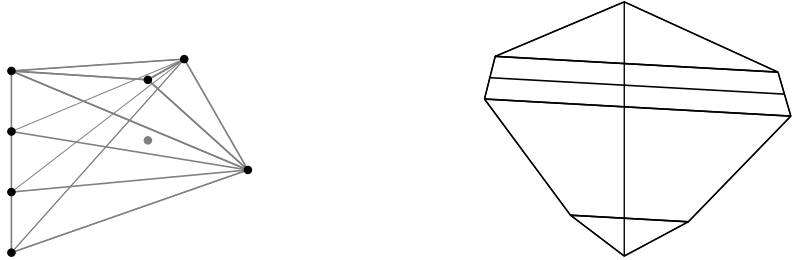
\begin{figure}[ht]
\tdplotsetmaincoords{35}{20}
\begin{tikzpicture}[line join=bevel,tdplot_main_coords, scale=1.4]

	\draw[gray, ultra thin] (1,0,0) -- (0,1,0) -- (-1,-1,-1) -- cycle;
	\draw[gray, semithick] (1,0,0) -- (0,1,0) -- (0,0,1) -- cycle;
	\draw[gray, semithick] (1,0,0) -- (0,0,1) -- (-1,-1,2) -- cycle;
	\draw[gray, ultra thin] (0,1,0) -- (-1,-1,-1) -- (-1,-1,2) -- cycle;
	\draw[gray, semithick] (0,1,0) -- (0,0,1) -- (-1,-1,2) -- cycle;
	\draw[gray, semithick] (1,0,0) -- (-1,-1,-1) -- (-1,-1,2) -- cycle;
	\draw[gray, semithick] (1,0,0) -- (-1,-1,0);
	\draw[gray, semithick] (1,0,0) -- (-1,-1,1);
	\draw[gray, ultra thin] (0,1,0) -- (-1,-1,0);
	\draw[gray, ultra thin] (0,1,0) -- (-1,-1,1);

	\filldraw[black] (1,0,0) circle (1pt);
	\filldraw[black] (0,1,0) circle (1pt);
	\filldraw[black] (0,0,1) circle (1pt);
	\filldraw[black] (-1,-1,-1) circle (1pt);
	\filldraw[black] (-1,-1,0) circle (1pt);
	\filldraw[black] (-1,-1,1) circle (1pt);
	\filldraw[black] (-1,-1,2) circle (1pt);
	\filldraw[gray] (0,0,0) circle (1pt);
\end{tikzpicture}
\hspace{2.4cm}
\tdplotsetmaincoords{50}{320}
\begin{tikzpicture}[line join=bevel,tdplot_main_coords,scale=1.1]

\draw[black, semithick] (-1,-1,-1) -- (-1,0,-1) -- (0,-1,-1) -- cycle;
\draw[black, ultra thin] (-1,-1,3) -- (-1,-1,-1) -- (0,-1,-1) -- (8/5,-1,-1/5) -- (7/5,-1,3/5) -- cycle;
\draw[black, ultra thin] (-1,-1,3) -- (-1,-1,-1) -- (-1,0,-1) -- (-1,8/5,-1/5) -- (-1,7/5,3/5) -- cycle;
\draw[black, semithick] (0,-1,-1) -- (-1,0,-1) -- (-1,8/5,-1/5) -- (8/5,-1,-1/5) -- cycle;
\draw[black, semithick] (-1,7/5,3/5) -- (-1,8/5,-1/5) -- (8/5,-1,-1/5) -- (7/5,-1,3/5) -- cycle;
\draw[black, semithick] (-1,-1,3) -- (-1,7/5,3/5) -- (7/5,-1,3/5) -- cycle;

\draw[black, semithick] (-1,3/2,1/5) -- (3/2,-1,1/5);

\end{tikzpicture}
\caption{The polyhedral subdivision and a schematic picture of $Y$}
\label{fig:34_decomposition}
\end{figure}

Figure~\ref{fig:34_decomposition} shows the unique fine mixed subdivision of the boundary of $P$.  This defines decomposition data for $P$.  Taking cones over the polyhedra in this decomposition gives a fan $\Sigma$ that defines a toric resolution~$Y$ of~$X_P$.  The variety $Y$ is smooth, with ten toric points and fifteen toric curves arranged as on the right-hand side of Figure~\ref{fig:34_decomposition}: the toric points are the vertices and the toric curves are the edges.  Note the two toric surfaces in $Y$ that map to the curve of singularities under the resolution $Y \to X_P$.

\begin{figure}[ht]
\tdplotsetmaincoords{35}{20}
\begin{tikzpicture}[line join=bevel,tdplot_main_coords,scale=1.4]

	\draw[gray, ultra thin] (1,0,0) -- (0,1,0) -- (-1,-1,-1) -- cycle;
	\draw[gray, semithick] (1,0,0) -- (0,1,0) -- (0,0,1) -- cycle;
	\draw[gray, semithick] (1,0,0) -- (0,0,1) -- (-1,-1,2) -- cycle;
	\draw[gray, ultra thin] (0,1,0) -- (-1,-1,-1) -- (-1,-1,2) -- cycle;
	\draw[gray, semithick] (0,1,0) -- (0,0,1) -- (-1,-1,2) -- cycle;
	\draw[gray, semithick] (1,0,0) -- (-1,-1,-1) -- (-1,-1,2) -- cycle;
	\draw[red, semithick] (1,0,0) -- (-1,-1,0);
	\draw[blue, semithick] (1,0,0) -- (-1,-1,1);
	\draw[gray, ultra thin] (0,1,0) -- (-1,-1,0);
	\draw[gray, ultra thin] (0,1,0) -- (-1,-1,1);

	\filldraw[black] (1,0,0) circle (1pt);
	\filldraw[black] (0,1,0) circle (1pt);
	\filldraw[black] (0,0,1) circle (1pt);
	\filldraw[black] (-1,-1,-1) circle (1pt);
	\filldraw[black] (-1,-1,0) circle (1pt);
	\filldraw[black] (-1,-1,1) circle (1pt);
	\filldraw[black] (-1,-1,2) circle (1pt);
	\filldraw[gray] (0,0,0) circle (1pt);
\end{tikzpicture}
\hspace{2.4cm}
\tdplotsetmaincoords{50}{320}
\begin{tikzpicture}[line join=bevel,tdplot_main_coords,scale=1.1]

\draw[black, semithick] (-1,-1,-1) -- (-1,0,-1) -- (0,-1,-1) -- cycle;
\draw[black, ultra thin] (-1,-1,3) -- (-1,-1,-1) -- (0,-1,-1) -- (8/5,-1,-1/5) -- (7/5,-1,3/5) -- cycle;
\draw[black, ultra thin] (-1,-1,3) -- (-1,-1,-1) -- (-1,0,-1) -- (-1,8/5,-1/5) -- (-1,7/5,3/5) -- cycle;
\draw[black, semithick] (0,-1,-1) -- (-1,0,-1) -- (-1,8/5,-1/5) -- (8/5,-1,-1/5) -- cycle;
\draw[black, semithick]  (-1,8/5,-1/5) -- (8/5,-1,-1/5) -- (7/5,-1,3/5) -- (-1,7/5,3/5);
\draw[black, semithick] (-1,-1,3) -- (-1,7/5,3/5) -- (7/5,-1,3/5) -- cycle;

\draw[black, semithick] (-1,3/2,1/5) -- (3/2,-1,1/5);
\draw[blue, semithick] (-1,8/5,-1/5) -- (-1,3/2,1/5);
\draw[red, semithick] (-1,7/5,3/5) -- (-1,3/2,1/5);

\end{tikzpicture}
\caption{Toric curves that generate the subspace of exceptional curves}
\label{fig:34_generators}
\end{figure}
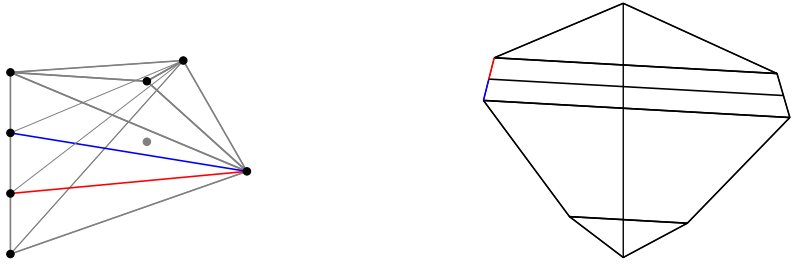

Applying Theorem~\ref{thm:toric_Betti} gives
\begin{align*}
	b_0(Y) &= 1  &  b_4(Y) &= 4   \\
	b_1(Y) &= 0  &  b_5(Y) &= 0  \\
	b_2(Y) &= 4 & b_6(Y) &= 1 \\
	b_3(Y) &= 0 &
\end{align*}
and since there are no quadrilaterals in the polyhedral subdivision, we find that the Betti numbers of $Y_t$ coincide with those of $Y$.  Computing the subspace $L \subset H_2(Y_t) \cong H_2(Y)$ of classes of exceptional curves, as in the cube example, we find that $l=2$ and that generators for $L$ are as shown in Figure~\ref{fig:34_generators}.  These generators are the fibers of the toric surfaces in $Y$ that resolve the curve of transverse $A_2$ singularities.

\noindent From \eqref{eq:N_e} and \eqref{eq:m} we find that (conjecturally) there are $m=3$ exceptional curves in total, and therefore:
\begin{align*}
	b_0(X_\eta) &= 1  &  b_4(X_\eta) &= 2   \\
	b_1(X_\eta) &= 0  &  b_5(X_\eta) &= 0  \\
	b_2(X_\eta) &= 2 & b_6(X_\eta) &= 1. \\
	b_3(X_\eta) &= 2 &
\end{align*}
These are the Betti numbers of the blow-up $X$ of $\PP^3$ in a plane cubic, which is consistent with the fact that the Minkowski polynomial
\[
f = \frac{(1+z)^3}{xyz} + x + y + z
\]
defined by our decomposition data for $P$ is a mirror to $X$.

\section{The Betti Numbers Depend Only on the Mirror Laurent Polynomial}
\label{sec:main_result}

We have seen that decomposition data for a three-dimensional reflexive polytope $P$ are, for each facet $F$ of $P$:
\begin{itemize}
\item[(A)] a choice of admissible Minkowski decomposition $F = F_1 + \cdots + F_k$;
\item[(B)] a choice of regular fine mixed subdivision of $F$ subordinate to (A);
\end{itemize}
and that decomposition data determine a Minkowski polynomial~$f$.  In this section we will show that the Betti numbers of the smoothing $X$ of $X_P$ determined by the decomposition data are independent of the choice (B).  This implies, in view of the discussion in \S\ref{sec:topology}, that the Betti numbers of $X$ depend on the decomposition data only via $f$.  We will prove:

\begin{theorem} \label{thm:main_result}
    Let $P$ be a three-dimensional reflexive polytope.  Consider two sets of decomposition data for $P$, where the choices (A) are the same in each case but the choices (B) are different.  Let $X$ denote the Corti--Hacking--Petracci smoothing of $X_P$ determined by the first set of decomposition data, with $Y \to X_P$ the toric partial resolution and $Y_t$ the smoothing of $Y$; let $X'$ denote the smoothing of $X_P$ determined by the second set of decomposition data, with $Y' \to X_P$ the toric partial resolution and $Y'_t$ the smoothing of~$Y'$.  Assume the Corti--Hacking--Petracci conjectures described on page~\pageref{fn:homology} hold. Then:
    \begin{enumerate}
        \item[(a)] $b_2(Y) = b_2(Y')$;
        \item[(b)] the Betti numbers of $Y_t$ and $Y_t'$ coincide;
        \item[(c)] the Betti numbers of $X$ and $X'$ coincide.
    \end{enumerate}
\end{theorem}

\noindent The key point is that the polyhedral decompositions of the boundary of $P$ that define $Y$ and $Y'$ differ by a sequence of the following two types of moves, or their inverses:
\begin{figure}[ht]
\centering
\begin{tikzpicture}

\draw[black, thick] (0,0) -- (1,0) -- (2,1) -- (1,2) -- (0,1) -- cycle;
\draw[black, thick] (0,1) -- (1,0);

\filldraw[black] (0,0) circle (1pt);
\filldraw[black] (0,1) circle (1pt);
\filldraw[black] (1,0) circle (1pt);
\filldraw[black] (1,1) circle (1pt);
\filldraw[black] (1,2) circle (1pt);
\filldraw[black] (2,1) circle (1pt);

\node at (3+1/2,1) {$\longrightarrow$};

\draw[black, thick] (5,0) -- (5+1,0) -- (5+2,1) -- (5+1,2) -- (5,1) -- cycle;
\draw[black, thick] (5,0) -- (5+1,1);
\draw[black, thick] (5+1,1) -- (5+1,2);
\draw[black, thick] (5+1,1) -- (5+2,1);

\filldraw[black] (5,0) circle (1pt);
\filldraw[black] (5,1) circle (1pt);
\filldraw[black] (5+1,0) circle (1pt);
\filldraw[black] (5+1,1) circle (1pt);
\filldraw[black] (5+1,2) circle (1pt);
\filldraw[black] (5+2,1) circle (1pt);

\end{tikzpicture}
\caption{A move of Type I}
\label{fig:Type_I}
\end{figure}
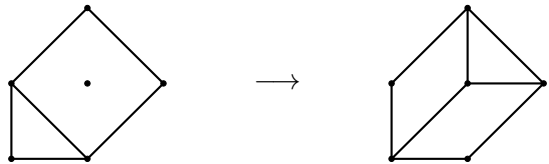

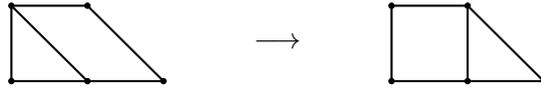
\begin{figure}[ht]
\centering
\begin{tikzpicture}

\draw[black, thick] (0,0) -- (1,0) -- (2,0) -- (1,1) -- (0,1) -- cycle;
\draw[black, thick] (0,1) -- (1,0);

\filldraw[black] (0,0) circle (1pt);
\filldraw[black] (0,1) circle (1pt);
\filldraw[black] (1,0) circle (1pt);
\filldraw[black] (1,1) circle (1pt);
\filldraw[black] (2,0) circle (1pt);

\node at (3+1/2,1/2) {$\longrightarrow$};

\draw[black, thick] (5,0) -- (5+1,0) -- (5+2,0) -- (5+1,1) -- (5+0,1) -- cycle;
\draw[black, thick] (5+1,1) -- (5+1,0);

\filldraw[black] (5,0) circle (1pt);
\filldraw[black] (5,1) circle (1pt);
\filldraw[black] (5+1,0) circle (1pt);
\filldraw[black] (5+1,1) circle (1pt);
\filldraw[black] (5+2,0) circle (1pt);

\end{tikzpicture}
\caption{A move of Type II}
\label{fig:Type_II}
\end{figure}

\noindent The vertices of the outer pentagon in Figure~\ref{fig:Type_I} are $(0,0)$, $(1,0)$, $(a+1,b)$, $(b,a+1)$, $(0,1)$, where $a$ and $b$ are positive coprime integers, and the interior lattice point pictured is at $(a,b)$; unless $a=b=1$ then there are other interior lattice points which are not pictured.  The precise values of $a$ and $b$ will not affect the analysis. The vertices of the outer quadrilateral in Figure~\ref{fig:Type_II} are at $(0,0)$, $(2,0)$, $(1,1)$, and $(0,1)$.

\subsection{Type~I Moves} \label{sec:TypeI}

Let us analyse how $b_2$ changes under a Type~I move. Suppose first that $Y_1$,~$Z$, and $Y_2$ are three-dimensional toric varieties defined by polyhedral decompositions of the boundary of $P$ that differ only in the following way\footnote{The co-ordinates of the vertices and interior lattice point pictured in Figure~\ref{fig:Y1ZY2} are as in Figure~\ref{fig:Type_I}.}:
\begin{figure}[ht]
\centering
\begin{tikzpicture}

\draw[black, thick] (0,0) -- (1,0) -- (2,1) -- (1,2) -- (0,1) -- cycle;
\draw[black, thick] (0,1) -- (1,0);

\filldraw[black] (0,0) circle (1pt);
\filldraw[black] (0,1) circle (1pt);
\filldraw[black] (1,0) circle (1pt);
\filldraw[gray] (1,1) circle (1pt);
\filldraw[black] (1,2) circle (1pt);
\filldraw[black] (2,1) circle (1pt);

\node at (2,1+3/4) {$Y_1$};

\draw[black, thick] (4,0) -- (4+1,0) -- (4+2,1) -- (4+1,2) -- (4,1) -- cycle;

\filldraw[black] (4,0) circle (1pt);
\filldraw[black] (4,1) circle (1pt);
\filldraw[black] (4+1,0) circle (1pt);
\filldraw[gray] (4+1,1) circle (1pt);
\filldraw[black] (4+1,2) circle (1pt);
\filldraw[black] (4+2,1) circle (1pt);

\node at (4+2,1+3/4) {$Z$};

\draw[black, thick] (8,0) -- (8+1,0) -- (8+2,1) -- (8+1,2) -- (8,1) -- cycle;
\draw[black, thick] (8,0) -- (8+1,1);
\draw[black, thick] (8+1,1) -- (8+1,2);
\draw[black, thick] (8+1,1) -- (8+2,1);

\filldraw[black] (8,0) circle (1pt);
\filldraw[black] (8,1) circle (1pt);
\filldraw[black] (8+1,0) circle (1pt);
\filldraw[black] (8+1,1) circle (1pt);
\filldraw[black] (8+1,2) circle (1pt);
\filldraw[black] (8+2,1) circle (1pt);

\node at (8+2,1+3/4) {$Y_2$};

\end{tikzpicture}
\caption{The fans for $Y_1$,~$Z$, and~$Y_2$ differ only at the cones over these polygons.}
\label{fig:Y1ZY2}
\end{figure}
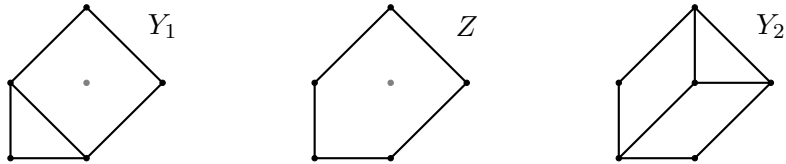

\begin{figure}[ht]
\centering
\begin{tikzpicture}

\draw[black, thick] (4,0) -- (4+1,0) -- (4+2,1) -- (4+1,2) -- (4,1) -- cycle;

\filldraw[black] (4,0) circle (1pt);
\filldraw[black] (4,1) circle (1pt);
\filldraw[black] (4+1,0) circle (1pt);
\filldraw[black] (4+1,2) circle (1pt);
\filldraw[black] (4+2,1) circle (1pt);

\node[label=right:{$a_1$}] at (4+2,1) {};
\node[label=above:{$a_2$}] at (4+1,2) {};
\node[label=left:{$a_3$}] at (4,1) {};
\node[label=below:{$a_4$}] at (4,0) {};
\node[label=below:{$a_5$}] at (4+1,0) {};

\end{tikzpicture}
\caption{The values of a piecewise-linear function on the fan for $Z_0$.}
\label{fig:line_bundle}
\end{figure}
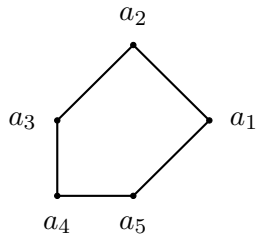

\noindent Then $Y_1$ and $Y_2$ are toric partial resolutions of $X_P$ which differ by a Type~I move, and there is a diagram
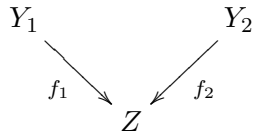
\[
\xymatrix{
    Y_1 \ar[dr]_{f_1} & & Y_2 \ar[dl]^{f_2} \\
    & Z
}
\]
Let $z \in Z$ be the torus-fixed point corresponding to the pentagon pictured, and write $Z^0 = Z \setminus \{z\}$.  The maps $f_1$ and $f_2$ induce isomorphisms $f_1^{-1}(Z^0) \to Z^0$ and $f_2^{-1}(Z^0) \to Z^0$, and the resulting inclusions $j_1 \colon Z^0 \to Y_1$ and $j_2 \colon Z^0 \to Y_2$ define a diagram
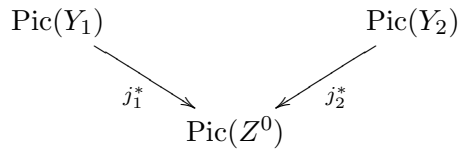
\[
\xymatrix{
    \Pic(Y_1) \ar[dr]_{j_1^*} & & \Pic(Y_2) \ar[dl]^{j_2^*} \\
    & \Pic(Z^0)
}
\]
We will identify $\Pic(Y_1)$ and $\Pic(Y_2)$ as subspaces of $\Pic(Z^0)$.  To give a line bundle on $Z^0$ is to give a piecewise-linear function on each maximal cone in the fan for $Z^0$, subject to the constraint that these piecewise-linear functions agree along faces.  Let us write the values of such a piecewise linear function at the vertices of the polyhedral decomposition that we are considering as in Figure~\ref{fig:line_bundle}.

Since the fan for $Z^0$ does not include the cone over the pentagon pictured in Figure~\ref{fig:line_bundle}, that cone does not impose any relation between the values $a_1,\ldots,a_5$.  (There may be other relations from the part of the fan not pictured, but these will be the same for $Z$, $Y_1$, and $Y_2$.)  This piecewise-linear function defines a line bundle on $Y_1$ if and only if it is piecewise-linear on the two cones pictured on the left-hand side of Figure~\ref{fig:Y1ZY2}, that is, if and only if $a_1+a_3 = a_2+a_5$.  So
\[
\Pic(Y_1) = \big\{ a_1 + a_3 = a_2 + a_5\} \subset \Pic(Z^0).
\]
Similarly, the piecewise-linear function on $Z_0$ defines a line bundle on $Y_2$ if and only if it is piecewise-linear on the two cones pictured on the right-hand side of Figure~\ref{fig:Y1ZY2}, that is, if and only if $a_1+a_4 - a_5= a_2+a_4 - a_3$. So
\[
\Pic(Y_2) = \big\{ a_1+a_4 - a_5= a_2+a_4 - a_3 \} \subset \Pic(Z^0).
\]
These are the same subspace of $\Pic(Z_0)$; therefore $\Pic(Y_1)$ and $\Pic(Y_2)$ are canonically isomorphic.  Since all maximal cones in the fan for $Y_1$ are full-dimensional, the Picard group of $Y_1$ is isomorphic to $H^2(Y_1)$; the same statement is true for $Z$ and for $Y_2$.  Thus $b_2$ is invariant under Type~I moves.

A result of Fulton and Sturmfels~\cite[Proposition~1.1]{FultonSturmfels1997} implies that $b_4(Y_i)$ is equal to three less than the number of vertices in the polyhedral decomposition that defines $Y_i$, and so $b_4(Y_2) = b_4(Y_1)+1$.  Furthermore $\chi(Y_2) = \chi(Y_1)+1$ -- here we used Theorem~\ref{thm:toric_Betti} -- and so $b_3(Y_1) = b_3(Y_2)$.  In summary:
\begin{align*}
    & b_0(Y_2) = b_0(Y_1) = 1 && b_4(Y_2) = b_4(Y_1) + 1\\
    & b_1(Y_2) = b_1(Y_1) = 0 && b_5(Y_2) = b_5(Y_1) = 0  \\
    & b_2(Y_2) = b_2(Y_1)  && b_6(Y_2) = b_6(Y_1) = 1 \\
    & b_3(Y_2) = b_3(Y_1)  &&
\end{align*}

\subsection{Type~II Moves}
An essentially identical argument shows that the Betti numbers of $Y$ are invariant under Type~II moves.  We are now in a position to prove Theorem~\ref{thm:main_result}.

\subsection{Proof of Theorem~\ref{thm:main_result}}

Since $Y$ and $Y'$ differ by a sequence of moves of Type~I and~II, and their inverses, we have that $b_2(Y') = b_2(Y)$.  This is part (a) of the Theorem. Furthermore $b_3(Y') = b_3(Y)$, and if the sequence of moves connecting $Y$ to $Y'$ contains $M$ moves of Type~I and $N$ moves of $(\text{Type I})^{-1}$ then $b_4(Y') = b_4(Y)+M-N$, and the numbers $k$ and $k'$ of quadrilaterals in the polyhedral decompositions defining $Y$ and $Y'$ satisfy $k' = k + M-N$.  The vanishing cycle analysis \eqref{eq:Y_to_Yt} now implies (b).

To prove (c), it suffices to show that the quantities $l$ and $m$ occurring in equation \eqref{eq:Yt_to_X} are the same for $Y$ and  $Y'$.  This is obvious for the number of nodes $m$: the conjectural formula \eqref{eq:m} for $m$ depends only on the sizes of the partitions of $\widetilde{\Gamma}_e$, which in turn depends on the choices of Minkowski decomposition (A) but not on the fine mixed subdivisions (B).  It remains to show that the dimension $l$ of  the subspace $L$ of $H_2(Y)$ spanned by the classes of exceptional curves in the smoothing $Y_t$ is the same as the dimension $l'$ of  the subspace $L'$ of $H_2(Y')$ spanned by the classes of exceptional curves in the smoothing $Y'_t$.  Let us return to the situation considered in \S\ref{sec:TypeI}, where $Y_1$ and $Y_2$ are three-dimensional toric varieties that differ by a Type~I move.  We showed there that $H^2(Y_1)$ and $H^2(Y_2)$ are isomorphic, via the inclusions
\[
\xymatrix{
    H^2(Y_1) \ar[dr]_{j_1^*} & & H^2(Y_2) \ar[dl]^{j_2^*} \\
    & H^2(Z^0).
}
\]
Dualising gives 
\[
\xymatrix{
    H_2(Y_1)  & & H_2(Y_2)  \\
    & H_2(Z^0) \ar[ul]^{{j_1}_*} \ar[ur]_{{j_2}_*}.
}
\]
and since the subspace of $H_2(Y_i)$ spanned by exceptional curves is pushed forward from $H_2(Z^0)$ via ${j_i}_*$ it follows that the dimension $l$ of this subspace is also invariant under Type~I moves.  Repeating this analysis for Type~II moves shows that $l=l'$, and proves Theorem~\ref{thm:main_result}.

\section{Systematic Analysis}
\label{sec:computations}

The computation of Betti numbers described in Section~\ref{sec:smoothing} can be automated.  The key ingredients are:
\begin{itemize}
\item algorithms for computing with lattice polyhedra and their duals.  There are several robust and well-tested implementations here, including those in Magma~\cite{Magma}, Sage~\cite{Sage}, and Polymake~\cite{Polymake2000}.
\item the Kreuzer--Skarke classification~\cite{KreuzerSkarke1998} of three-dimensional reflexive polytopes.
\item Altmann's determination~\cite{Altmann1997} of all Minkowski summands of a given polytope.
\item the computation of fine mixed subdivisions (Definition~\ref{def:subdivision}), that is, the determination of all regular triangulations of a Cayley polytope.  For this we use J\"org Rambau's TOPCOM package~\cite{TOPCOM}.
\item an HPC cluster.  Some of the computations involved are quite challenging.
\end{itemize}
\label{sec:systematic}
Full source code for these computations, written in Magma, can be found at~\cite{this_repo}.  This relies in an essential way on code from the Fanosearch project~\cite{fanosearch_repo}.

There are 4319 three-dimensional reflexive polytopes, which in total admit more than a billion decomposition data.  These decomposition data give rise to 3857 distinct Minkowski polynomials\footnote{The number of Minkowski polynomials here differs slightly from the count in~\cite{AkhtarCoatesGalkinKasprzyk2012}, because there the authors required Minkowski decompositions of facets to satisfy an additional lattice condition (ibid., Definition~7) and here we do not regard $\GL(3,\ZZ)$-equivalent Minkowski polynomials as the same.}, which together give mirrors to the 98 three-dimensional Fano manifolds with very ample\footnote{The seven three-dimensional Fano manifolds without very ample canonical bundle are not expected to admit Laurent polynomial mirrors with reflexive Newton polytopes, and so fall outside the range of the Corti--Hacking--Petracci construction.} anticanonical bundle.  We analysed several million decomposition data, including at least one decomposition for each of the 3857 Minkowski polynomials.  In each case we found that
\begin{equation} \label{eq:statement}
    \tag{$\ast$} \begin{minipage}{0.9\textwidth} the Betti numbers of the smoothing $X_\eta$ determined by the decomposition data depend only on the mirror Fano manifold $X$, and coincide\footnotemark[5] with those of $X$.\end{minipage}
\end{equation}   
This provides evidence for the conjectural picture described in the introduction:\footnotetext[5]{Betti numbers for three-dimensional Fano manifolds can be found in~\cite{IskovskikhProkhorov99}.} that if a Fano manifold $X$ corresponds under Mirror Symmetry to
a Laurent polynomial $f$ then there is a degeneration $\cX \to \Delta$ with general fiber $X$ and special fiber the toric variety defined by the
spanning fan of the Newton polytope of $f$. It also provides evidence for the conjectures of Corti--Hacking--Petracci described in
Section~\ref{sec:smoothing}, on the number and homology class of the exceptional curves in their resolution $\pi_t \colon Y_y \to X_t$. If these
conjectures are correct then, in view of Theorem~\ref{thm:main_result}, 3857 of these calculations give a computer-assisted rigorous proof of \eqref{eq:statement}.

\section*{Acknowledgements}

This project has received funding from the European Research Council (ERC) under the European Union’s Horizon 2020 research and innovation programme (grant agreement No.~682603), and from the EPSRC Programme Grant EP/N03189X/1, \emph{Classification, Computation, and Construction: New Methods in Geometry}.  We thank Paul Hacking, Andrea Petracci, and Alexander Kasprzyk for many extremely useful conversations, and thank Andy Thomas, Matt Harvey, and the Imperial College Research Computing Service team for invaluable technical assistance.

\bibliographystyle{amsplain}
\bibliography{bibliography}

\end{document}